\documentclass[11pt]{amsart}
\usepackage{amssymb,amscd,amsmath,amsthm}
\usepackage[all]{xy}

\theoremstyle{plain}
\newtheorem{theorem}{Theorem}[section]
\newtheorem{proposition}[theorem]{Proposition}

\newtheorem{lemma}[theorem]{Lemma}

\theoremstyle{definition}
\newtheorem{definition}[theorem]{Definition}

\newtheorem{remark}[theorem]{Remark}

\numberwithin{equation}{section}

\newcommand{\CC}{{\mathbb{C}}}

\newcommand{\RR}{{\mathbb{R}}}
\newcommand{\SSS}{{\mathbb{S}}}

\newcommand{\ZZ}{{\mathbb{Z}}}

\newcommand{\id}{{\operatorname{id}}}

\newcommand{\SO}{{\operatorname{SO}}}

\newcommand{\Diff}{{\operatorname{Diff}}}
\newcommand{\sm}{{\smallsetminus}}
\newcommand{\eps}{\varepsilon}
\newcommand{\e}{\mathrm{e}}

\begin{document}

\pagestyle{plain}

\title{The orientation-preserving diffeomorphism group of $\SSS^2$ deforms to $\SO(3)$ smoothly}

\author{Jiayong Li}
\email{jiayong.li@utoronto.ca}
\address{Department of Mathematics, University of Toronto, Room 6290, 40 St. George Street, Toronto, Ontario, Canada M5S 2E4}

\author{Jordan Alan Watts}
\email{jwatts@math.toronto.edu}
\address{Department of Mathematics, University of Toronto, Room 6290, 40 St. George Street, Toronto, Ontario, Canada M5S 2E4}

\begin{abstract}
Smale proved that the orientation-preserving diffeomorphism group of $\SSS^2$ has a \emph{continuous} strong deformation retraction to $\SO(3)$. In this paper, we construct such a strong deformation retraction which is \emph{diffeologically smooth}.
\end{abstract}

\maketitle

\section{Introduction}
In Smale's 1959 paper ``Diffeomorphisms of the 2-Sphere" (\cite{smale}), he shows that there is a continuous strong deformation retraction from the orientation-preserving $C^\infty$ diffeomorphism group of $\SSS^2$ to the rotation group $\SO(3)$. The topology of the former is the $C^k$ topology. In this paper, we construct such a strong deformation retraction which is diffeologically smooth. We follow the general idea of \cite{smale}, but to achieve smoothness, some of the steps we use are completely different from those of \cite{smale}. The most notable differences are explained in Remark \ref{comparison1}, \ref{comparison2}, and Remark \ref{comparison3}. We note that there is a different proof of Smale's result in \cite{earle}, but the homotopy is not shown to be smooth.

\

We start by defining the notion of diffeological smoothness in three special cases which are directly applicable to this paper. Note that diffeology can be defined in a much more general context and we refer the readers to \cite{iglesias}.

\begin{definition} \label{plot}
Let $U$ be an arbitrary open set in a Euclidean space of arbitrary dimension.
\begin{itemize}
\item
Suppose $\Lambda$ is a manifold with corners. A map $P: U \to \Lambda$ is a \textbf{plot} if $P$ is $C^\infty$.
\item
Suppose $X$ and $Y$ are manifolds with corners, and $\Lambda \subset C^\infty (X,Y)$. Denote by $\mathrm{ev}$ the evaluation map $\Lambda \times X \to Y$ given by $(f,x) \mapsto f(x)$. A map $P: U \to \Lambda$ is a \textbf{plot}, if the map from $U \times X$ to $Y$ given by $(s,x) \mapsto \mathrm{ev}(P(s),x)$ is $C^\infty$.
\item
Suppose $\Lambda$ is given by a product $\Lambda_1 \times \cdots \times \Lambda_n$, where each $\Lambda_i$ is either one of the above two cases. Denote by $\pi_i$ the projection onto each factor. A map $P: U \to \Lambda$ is a \textbf{plot} if each map $\pi_i \circ P: U \to \Lambda_i$ is a plot.
\end{itemize}
\end{definition}

In the second case, one can think of a plot $P$ of $\Lambda$ as a smooth family of maps $\{f_s\}_{s \in U} \subset \Lambda$ by letting $f_s = P(s)$.

\begin{definition}
Let $\Lambda$ be one of the three cases in Definition \ref{plot}. Then $\Lambda$ equipped with the collection of all the plots is called a \textbf{diffeological space}. The second kind is referred to as the standard functional diffeology, and the third kind the product diffeology.
\end{definition}

Let $\Lambda$ and $\Gamma$ be two diffeological spaces. We define the notion of diffeological smoothness of a map from $\Lambda$ to $\Gamma$ as follows.
\begin{definition}
A map $\varphi: \Lambda \to \Gamma$ is \textbf{diffeologically smooth} if for every plot $P: U \to \Lambda$, the map $\varphi \circ P: U \to \Gamma$ is a plot. If in addition, $\varphi$ has a diffeological smooth inverse, then it is a \textbf{diffeomorphism}.
\end{definition}

\begin{remark}
Let $\Lambda$ and $\Gamma$ be manifolds with corners. Then $\varphi: \Lambda \to \Gamma$ is a diffeologically smooth map if and only if $\varphi$ is a smooth map between manifolds with corners. We refer the readers to \cite{joyce} for the definition of smooth maps between manifolds with corners.

This fact is non-trivial. The proof of ``IV.13 Smooth real maps from half-spaces" in \cite{iglesias} contains the main ingredients for showing the equivalence of two notions of smoothness for manifolds with boundary. One can use Theorem 1 of \cite{schwarz} (take the group $G$ to be $(\ZZ/2)^n$) to modify this proof to show the equivalence for manifolds with corners.

Therefore, we see that diffeological smoothness is a generalization of the usual notion of smoothness. Throughout the rest of the paper, by a ``smooth map" we always mean a diffeologically smooth map.
\end{remark}

\

Now we state the main result of this paper. Let us denote
$$\Omega := \text{ the orientation-preserving $C^\infty$ diffeomorphism group of $\SSS^2$.}$$
\begin{theorem}[main theorem]\label{main}
There is a smooth strong deformation retraction $P:I\times\Omega\to\Omega$ to $\SO(3)$ that is equivariant under the left action of $\SO(3)$. More precisely, for each $(t,f) \in I \times \Omega$ and $A \in \SO(3)$,
\begin{enumerate}
\item
$P_0(f) = f$,
\item
$P_1(f) \in \SO(3)$,
\item
$P_t(A) = A$,
\item
$P_t(A \circ f) = A \circ P_t(f)$.
\end{enumerate}
\end{theorem}

Now we outline the construction of the above deformation retraction $P$. Let $x_0$ be the South Pole of $\SSS^2 \subset \RR^3$, and $\e_1 = (1,0,0)$ and $\e_2 = (0,1,0)$ the basis vectors of the tangent space $T_{x_0}\SSS^2$. Denote by $\Omega_1$ the following subset of $\Omega$.
$$\Omega_1 := \{ f \in \Omega: f(x_0) = x_0 \text{ and } df|_{x_0} = \id_{T_{x_0}\mathbb{S}^2} \}.$$
\begin{lemma} \label{embed}
The map $i: \SO(3) \times \Omega_1 \to \Omega$ given by $(A,f) \mapsto A^{-1} \circ f$ is a diffeomorphism with its image
$$i(\SO(3) \times \Omega_1) = \{ f \in \Omega: df|_{x_0}\e_1 \text{ and } df|_{x_0}\e_2 \text{ are orthonormal} \}.$$
\end{lemma}

In the following theorem, we homotope $\Omega$ to $i(\SO(3) \times \Omega_1)$.
\begin{theorem} \label{Omega1}
There is a smooth homotopy $Q: I \times \Omega \rightarrow \Omega$ to $i(\SO(3) \times \Omega_1)$ that fixes $\SO(3)$ and is equivariant under the left action of $\SO(3)$. More precisely, for each $(t,f) \in I \times \Omega$ and $A \in \SO(3)$,
\begin{enumerate}
\item
$Q_0(f) = f$,
\item
$Q_1(f) \in i(\SO(3) \times \Omega_1)$,
\item
$Q_t(A) = A$,
\item
$Q_t(A \circ f) = A \circ Q_t(f)$.
\end{enumerate}
\end{theorem}

In the following theorem, we homotope $\Omega_1$ to $\{\id_{\mathbb{S}^2}\}$.
\begin{theorem} \label{Id}
There is a smooth strong deformation retraction $R: I \times \Omega_1 \rightarrow \Omega_1$ to $\{ \id_{\mathbb{S}^2} \}$. More precisely, for each $(t,f) \in I \times \Omega_1$,
\begin{enumerate}
\item
$R_0(f) = f$,
\item
$R_1(f) = \id_{\SSS^2}$,
\item
$R_t(\id_{\SSS^2}) = \id_{\SSS^2}$.
\end{enumerate}
\end{theorem}

Smoothly concatenating homotopies $Q$ and $R$ gives the desired deformation retraction $P$.

\

When proving Theorem \ref{Id}, we need a proposition about the diffeomorphisms of the square $[-1,1]^2$. Let $\mathcal F$ be the space of those orientation-preserving diffeomorphisms of the square $[-1,1]^2$ such that for each $f \in \mathcal F$, there exists a neighborhood of the boundary $\partial ([-1,1]^2)$ on which $f$ is the identity map.
\begin{proposition} \label{F}
There is a smooth strong deformation retraction $F: I \times \mathcal F \to \mathcal F$ to $\{\id_{I^2}\}$. More precisely, for each $(t,f) \in I \times \mathcal F$,
\begin{enumerate}
\item
$F_0(f) = f$,
\item
$F_1(f) = \id_{[-1,1]^2}$,
\item
$F_t(\id_{[-1,1]^2}) = \id_{[-1,1]^2}$.
\end{enumerate}
\end{proposition}

This paper is organized as follows. In Section 2, we construct the homotopy of Theorem \ref{Omega1}, and prove Theorem \ref{Id} while assuming Proposition \ref{F}. Lastly we prove the main theorem. In Section 3, we prove Proposition \ref{F}.

\

\textbf{Acknowledgement.} We would like to express our deepest gratitude to Yael Karshon for her time and patience involved in supervising this project. We would also like to thank Katrin Wehrheim for a helpful suggestion in the proof of Lemma \ref{epsilon}.

\

\section{Construction of homotopies}
We start by surveying some commonly used properties regarding diffelogical smoothness. One convenience of working with diffeology on $C^\infty$ function spaces is that verifying diffeological smoothness often reduces to checking the usual smoothness in finite dimensions. The proofs of these properties are left as exercises.
\begin{remark} \label{property}

\

\begin{enumerate}
\item
Let $\Lambda$, $\Gamma$, and $\Sigma$ be diffeological spaces. If $\varphi: \Lambda \to \Gamma$ and $\psi: \Gamma \to \Sigma$ are both diffeologically smooth, then $\psi \circ \varphi$ is diffeologically smooth.
\item
Let $X$, $Y$ and $Z$ be manifolds with corners. Let $\Lambda$ be a subset of $C^\infty (X,Y)$, and $\Gamma$ a subset of $C^\infty (Y,Z)$. Then the map from $\Lambda \times \Gamma$ to $C^\infty (X,Z)$ given by $(f,g) \mapsto g \circ f$ is diffeologically smooth.
\item
Let $X$ be a manifold with corners, and denote by $\Diff(X)$ the diffeomorphism group of $X$. Let $\Lambda$ be a subset of $\Diff(X)$. Then the map from $\Lambda$ to $\Diff(X)$ given by $f \mapsto f^{-1}$ is diffeologically smooth. (This fact can be proved by using the usual implicit function theorem.)
\item
Let $X$ and $Y$ be manifolds with corners, and $\Lambda$ a subset of $C^\infty (X,Y)$. Then the map from $\Lambda$ to $C^\infty (TX,TY)$ given by $f \mapsto Tf$ is diffeologically smooth.
\end{enumerate}
\end{remark}

Now we prove Lemma \ref{embed}. Recall that $x_0$ is the South Pole of $\SSS^2$, and $\e_1 = (1,0,0)$ and $\e_2 = (0,1,0)$ are the basis vectors of the tangent space $T_{x_0}\SSS^2$.
\begin{proof}[Proof of Lemma \ref{embed}]
Denote by $\tilde \Omega$ the following set $$\tilde \Omega := \{ f \in \Omega: df|_{x_0}\e_1 \text{ and } df|_{x_0}\e_2 \text{ are orthonormal} \}.$$ It follows from property (2) of Remark \ref{property} that the map $i: \SO(3) \times \Omega_1 \to \Omega$ given by $(A,f) \mapsto A^{-1} \circ f$ is smooth. Moreover, tangent vectors $d(A^{-1} \circ f)|_{x_0}\e_1$ and $d(A^{-1} \circ f)|_{x_0}\e_2$ are clearly orthonormal. Thus $i(\SO(3) \times \Omega_1) \subset \tilde \Omega$.

Define the map $\alpha: \tilde \Omega \to \SO(3)$ as follows. For $g \in \tilde \Omega$, let $\alpha(g)$ be the element in $\SO(3)$ that sends the ordered basis $\{g(x_0),dg(x_0)\e_1,dg(x_0)\e_2\}$ to the ordered basis $\{x_0,\e_1,\e_2\}$. Note that $\alpha(g)$ can be written as a matrix involving the partial derivatives of $g$, and is easily seen to be smooth with respect to $g$. It is easy to see that $\alpha(g) \circ g$ is an element of $\Omega_1$.

It is straightforward to check that the map from $\tilde \Omega$ to $\SO(3) \times \Omega_1$ given by $g \mapsto (\alpha(g), \alpha(g) \circ g)$ is the smooth inverse of the map $i$.
\end{proof}

\

Now we construct the smooth homotopy from $\Omega$ to
$$i(\SO(3) \times \Omega_1) = \{ f \in \Omega: df|_{x_0}\e_1 \text{ and } df|_{x_0}\e_2 \text{ are orthonormal} \}.$$
The construction involves some computation in local coordinates. Let $p: \mathbb{S}^2 \sm \{-x_0\} \rightarrow \mathbb{R}^2$ be the stereographic projection from the North Pole. More precisely, $p(x_1,x_2,x_3) = (\frac{2x_1}{1-x_3}, \frac{2x_2}{1-x_3})$. It is convenient to denote by $B(\delta)$ the open ball in $\RR^2$ around $0$ with radius $\delta$, that is,
$$B(\delta):= \{ y\in \RR^2: |y|<\delta \}.$$

\begin{proof}[Proof of Theorem \ref{Omega1}]
Given $f\in\Omega$, define tangent vectors $u$ and $v$ in $T_{f(x_0)}\SSS^2$ by $u=df|_{x_0}\e_1$ and $v=df|_{x_0}\e_2$. Let $u_1=\frac{u}{||u||}$ be the normalization of $u$, and $u_2$ the vector which makes $\{ u_1, u_2 \}$ a positively oriented orthonormal basis in $T_{f(x_0)}\SSS^2$.

Define the map $\alpha: \Omega \to \SO(3)$ as follows. For $f \in \Omega$, let $\alpha(f)$ be the unique element in $\SO(3)$ that sends the ordered basis $\{f(x_0),u_1,u_2\}$ to the ordered basis $\{x_0,\e_1,\e_2\}$. Define the map $p_f: \SSS^2 \sm \{ -f(x_0) \} \to \RR^2$ by
$$p_f = p \circ \alpha(f).$$
The map $p_f$ is a coordinate chart with $p_f(f(x_0)) = 0$. The expression of the vector $u_1$ in this coordinate chart is $(1,0)$. More precisely,
$$dp_f|_{f(x_0)}u_1 = dp|_{x_0}\alpha(f)u_1
= \begin{pmatrix}
1 & 0 & 0 \\
0 & 1 & 0
\end{pmatrix}
\begin{pmatrix}
1 \\ 0 \\ 0
\end{pmatrix}
= \begin{pmatrix} 1 \\ 0 \end{pmatrix}.$$
Similarly, $dp_f|_{f(x_0)}u_2 = (0,1)$. By abuse of notation, we identify a vector $w \in T_{f(x_0)}\SSS^2$ with its expression in the coordinate chart $p_f$. Therefore,
$$
u_1 =
\begin{pmatrix}
1 \\ 0
\end{pmatrix}, \
u_2 =
\begin{pmatrix}
0 \\ 1
\end{pmatrix}, \
u =
\begin{pmatrix}
a \\ 0
\end{pmatrix}, \
v =
\begin{pmatrix}
b \\ c
\end{pmatrix},
$$
where $a$, $b$, and $c$ depend smoothly on $f$. By definition, $a$ is positive. Since $f$ is orientation-preserving, $c$ is also positive. Define the matrix $g_{f,1}$ by $$g_{f,1} =
\begin{pmatrix}
a & b \\
0 & c
\end{pmatrix}^{-1}.
$$

By construction, $g_{f,1}$ sends $u$ to $u_1$ and $v$ to $u_2$. Moreover, $g_{f,1}$ is upper triangular and with positive diagonal entries, and it depends smoothly on $f$. For each $t\in[0,1]$, define the matrix $g_{f,t}$ by the following linear interoplation.
$$g_{f,t} = (1-t)\id + tg_{f,1}.$$
It is easy to see that the matrix $g_{f,t}$ is an orientation-preserving isomorphism and depends smoothly on $(t,f)$. Thus the family of diffeomorphisms $g_{f,t}:\RR^2 \to \RR^2$ determines a time-dependent vector field $X_{f,t}$ on $\RR^2$. More precisely, the vector field $X_{f,t}$ is defined by
$$\frac{d}{dt}g_{f,t} = X_{f,t} \circ g_{f,t}.$$
Let $\rho: \RR^2 \to [0,1]$ be a smooth function such that $\rho(y)$ is $1$ for $y \in B(1)$ and $0$ for $y \notin B(2)$. Define the vector field $\tilde X_{f,t}$ by $\tilde X_{f,t} = \rho X_{f,t}$. This vector field determines a family of orientation-preserving diffeomorphisms $\tilde g_{f,t}$, which is given by the solution to the following ordinary differential equation.
$$\frac{d}{dt}\tilde g_{f,t} = \tilde X_{f,t} \circ \tilde g_{f,t}, \text{with the initial condition }\tilde g_{f,0}=\id.$$
It follows from the standard regularity theory of ordinary differential equations that $\tilde g_{f,t}$ is smooth with respect to $(t,f)$. It is easy to see that the restriction of $\tilde g_{f,t}$ to the complement of $B(2)$ is the identity map. Moreover, for each $t$, the restriction of $\tilde g_{f,t}$ to the set $\tilde B_f := \bigcap_{t\in[0,1]}\tilde g_{f,t}^{-1}(B(1))$ agrees with the restriction of $g_{f,t}$. One can show that $\tilde B_f$ contains an open neighborhood of the origin by using the compactness of the interval $[0,1]$. In particular, this shows that $d\tilde g_{f,1}|_0 = g_{f,1}$. Define $\Psi_{f,t}: \SSS^2 \sm \{-f(x_0) \} \to \SSS^2 \sm \{-f(x_0) \}$ by
$$\Psi_{f,t} = p_f^{-1} \circ \tilde g_{f,t} \circ p_{f}.$$
Since $\tilde g_{f,t}$ is a diffeomorphism which agrees with the identity map on the complement of $B(2)$, defining $\Psi_{f,t}(-f(x_0)) = -f(x_0)$ clearly yields an orientation-preserving diffeomorphism of $\SSS^2$. It is clear that $\Psi_{f,0} = \id_{\SSS^2}$. Furthermore, it follows from $d\tilde g_{f,1}|_0 = g_{f,1}$ that the differential $d\Psi_{f,1}|_{f(x_0)}$ sends $u$ to $u_1$ and $v$ to $u_2$.

Now define the deformation retraction $Q: I \times \Omega \to \Omega$ by
$$Q_t(f) = \Psi_{f,t} \circ f.$$
The map $Q$ is indeed well-defined since each $\Psi_{f,t} \circ f \in \Omega$. The smoothness of $Q$ follows from the fact that $\Psi_{f,t}$ is smooth with respect to $(t,f)$. Since $\Psi_{f,0} = \id_{\SSS^2}$, we have $Q_0(f) = f$. Moreover, the differential $d(\Psi_{f,1} \circ f)|_{x_0}$ sends $\e_1$ to $u_1$ and $\e_2$ to $u_2$. Hence $Q_1(f) \in i(\SO(3) \times \Omega_1)$.

For $f \in \Omega$ such that $df|_{x_0}\e_1$ and $df|_{x_0}\e_2$ are orthonormal, the corresponding $g_{f,t} = \id$, which implies that $\Psi_{f,t} = \id_{\SSS^2}$. This shows that the homotopy $Q$ fixes $i(\SO(3) \times \Omega_1)$. In particular, $Q$ fixes $\SO(3)$.

For $A \in \SO(3)$, to show that $Q_t(A \circ f) = A \circ Q_t(f)$, we first observe that $\alpha(A \circ f) = \alpha(f) \circ A^{-1}$. Hence $p_{A\circ f} = p_f \circ A^{-1}$, which implies
$$Q_t(A \circ f) = A \circ p_f^{-1} \circ \tilde g_{A\circ f,t} \circ p_f \circ f.$$
By using the definition of local charts $p_{A\circ f}$ and $p_f$, it is easy to check that the matrix $g_{A \circ f,1}$ is equal to the matrix $g_{f,1}$. Then it follows that $g_{A \circ f,t} = g_{f,t}$, which in turn implies that $\tilde g_{A \circ f,t} = \tilde g_{f,t}$. This finishes proving Theorem \ref{Omega1}.
\end{proof}

\begin{remark} \label{alpha}
Note that the map $\alpha$ defined in the proof of Lemma \ref{embed} is the restriction of $\alpha$ defined in the proof of Theorem \ref{Omega1}.
\end{remark}

\begin{remark} \label{comparison1}
In the above construction, $\tilde g_{f,t}$ is the diffeomorphism of $\RR^2$ which agrees with the linear map $g_{f,t}$ around the origin, and is the identity map outside a ball. This map is inspired by the diffeomorphism $G_{\nu}(f)$ on page 626 of \cite{smale}. However, this diffeomorphism is not entirely accurate. Using our notation, the diffeomorphism $G_{t}(f): \RR^2 \to \RR^2$ is defined by
$$G_{t}(f)(x) = \gamma(f,x) \ g_{f,t}(x) + (1-\gamma(f,x))x,$$
where $\gamma(f,x)$ is a function which is $1$ on a ball $B(\frac{1}{2}\eps(f))$, and $0$ outside a larger ball $B(\eps(f))$. It is clear that the map $G_{t}(f)$ agrees with $g_{f,t}$ on the smaller ball $B(\frac{1}{2}\eps(f))$, and is the identity map outside $B(\eps(f))$. However, part of the image $g_{f,t}(B(\frac{1}{2}\eps(f)))$ could lie outside of the larger ball $B(\eps(f))$, which makes $G_{t}(f)$ not injective. Even though such problem can be fixed by choosing a different $\gamma$, it is still not clear that $\gamma$ can be chosen so that $G_{t}(f)$ is a diffeomorphism. Thus instead of trying to construct such a function $\gamma$, we use the vector field approach to construct the diffeomorphism $\tilde g_{f,t}$.
\end{remark}

\

Now we prove Theorem \ref{Id}, that is, there is a smooth strong deformation retraction from the set $$\Omega_1 := \{ f \in \Omega: f(x_0) = x_0 \text{ and } df|_{x_0} = \id_{T_{x_0}\mathbb{S}^2} \}$$
to $\{ \id_{\mathbb{S}^2} \}$. We construct this deformation retraction in the following two steps.
\begin{enumerate}
\item
Homotope each diffeomorphism $f$ in $\Omega_1$ to a diffeomorphism which is the identity map on a neighborhood of the South Pole $x_0$.
\item
Homotope the diffeomorphism $f$ to the identity map on the complement of the neighborhood.
\end{enumerate}

We start by proving a technical lemma about choosing such a neighborhood.

Let $p: \mathbb{S}^2 \sm \{-x_0\} \rightarrow \mathbb{R}^2$ be the stereographic projection from the North Pole, and $(y_1,y_2)$ the coordinate variables of $\mathbb{R}^2$. Using $p$ as a local coordinate chart, we denote the local expression of the map $f \in \Omega_1$ by $\bar f$, which is given by
$$\bar f := p \circ f \circ p^{-1}.$$
Note that $\bar f$ is not necessarily well-defined on all of $\mathbb{R}^2$; if $f^{-1}(-x_0)$ is not equal to $-x_0$, then $\bar f$ is not defined at the point $p(f^{-1}(-x_0))$. However, for $f \in \Omega_1$, the local expression $\bar f$ is always defined at $0$, and we have $\bar f(0) = 0$. In addition, $d\bar f|_0 = \id$. Thus we can choose a neighborhood of $0$ on which $\bar f$ is well-defined and the operator $d\bar f - \id$ is small in the operator norm. The following technical lemma says that such choice of neighborhood can be made smoothly with respect to $f$.

\begin{lemma} \label{epsilon}
There exists a smooth function $\eps: \Omega_1 \rightarrow (0, \frac{1}{2}]$ such that \begin{itemize}
\item $\bar f$ is well-defined on the open ball $B(2\eps(f))$,
\item for $y \in B(\eps(f))$, we have $\|d\bar f|_{y} - \id \| < \frac{1}{4}$,
\end{itemize}
where $\| \cdot \|$ denotes the operator norm induced by the Euclidean norm on $\mathbb{R}^2$.
\end{lemma}

\begin{remark} \label{comparison2}
Later in the construction, the function $\eps$ is involved in the deformation retraction from $\Omega_1$ to $\{\id_{\mathbb{S}^2}\}$. Since the deformation retraction needs to depend smoothly on $f \in \Omega_1$, the function $\eps$ needs to be smooth. Note that the topological construction of $\eps$ on the top of page 625 of \cite{smale} is continuous but not smooth in general. As a result, our method of constructing $\eps$ is completely different. We use the Sobolev inequality to construct the smooth function $\eps$.
\end{remark}

\begin{proof}[Proof of Lemma \ref{epsilon}]
Let $h$ be the height function $h: \mathbb{S}^2 \rightarrow [-1,1]$ given by $h(x_1,x_2,x_3)=x_3$. Define $h_f: \overline{B(1)} \rightarrow [-1,1]$ by $h_f = h \circ f \circ p^{-1}|_{\overline{B(1)}}$. It is easy to see that $f \mapsto h_f$ is smooth. For $f \in \Omega_1$, the definition of $\Omega_1$ implies that $h_f(0) = h(x_0)=-1$. We also observe that the value $\bar f(y)$ is well-defined if and only if $h_f(y) < 1$. We need to use this criterion later in the proof.

By the Sobolev inequality (see Theorem 6, page 270 of \cite{evans}), there is a universal constant $c$ such that
\begin{equation} \label{sobolev}
\|u\|_{C^1\left(\overline{B(1)}\right)} \leq c \|u\|_{H^3(B(1))},
\end{equation}
for every smooth function $u$ on the closed unit ball $\overline{B(1)}$.

Here $\|u\|_{C^1\left(\overline{B(1)}\right)}$ is the $C^1$ norm of $u$ on $\overline{B(1)}$, defined by
$$\|u\|_{C^1\left(\overline{B(1)}\right)} := \sup_{\overline{B(1)}}|u| + \sup_{\overline{B(1)}}|\partial_{y_1}u| + \sup_{\overline{B(1)}} |\partial_{y_2}u|.$$

$\|u\|_{H^3(B(1))}$ is the Sobolev-3 norm of $u$ on $B(1)$, defined by $$\|u\|_{H^3(B(1))} := (\sum_{|\alpha| \leq 3} \int_{B(1)} |\partial^{\alpha}u|^2)^{\frac{1}{2}},$$ where $\alpha = (\alpha_1,\alpha_2)$ with $\alpha_1$ and $\alpha_2$ non-negative integers, $|\alpha| = \alpha_1+\alpha_2$, and $\partial^\alpha = \partial^{\alpha_1}_{y_1} \partial^{\alpha_2}_{y_2}$.

For $f\in \Omega_1$, define $$\eps_1(f) = \frac{1}{c}\left[\frac{1}{c^2} + \sum_{|\alpha| \leq 3} \int_{B(1)} |\partial^{\alpha}h_f|^2\right]^{-\frac{1}{2}}.$$

We claim that $\eps_1$ is a smooth function from $\Omega_1$ to $(0,1]$, and $\bar f$ is well-defined on $B(\eps_1(f))$. It is eay to check that $\eps_1(f) > 0$ and $\eps_1(f) \leq \frac{1}{c}c = 1$. The smoothness of $\eps_1$ follows from the smoothness of $h_f$ with respect to $f$ and the smoothness property of integration. Moreover, for $y \in B(\eps_1(f)) \subset B(1)$, the mean value theorem implies that $$|h_f(y) - h_f(0)| \leq \sup_{B(1)} |d(h_f)| \cdot |y|.$$ Thus for $y$ such that $|y| < \eps_1(f)$,
\begin{align*}
h_f(y) &\leq h_f(0) + \sup_{B(1)} |d(h_f)| \cdot |y| \\
&< -1 + \frac{\sup_{B(1)} |d(h_f)|}{c\left[\frac{1}{c^2} + \sum_{|\alpha| \leq 3} \int_{B(1)} |\partial^{\alpha}h_f|^2\right]^{\frac{1}{2}}} \\
&\leq -1 + \frac{\|h_f\|_{C^1\left(\overline{B(1)}\right)}}{c\|h_f\|_{H^3(B(1))}} \leq 0 < 1,
\end{align*}
where the second inequality is by definition of $\eps_1$ and the last line is by the Sobolev inequality (\ref{sobolev}). Therefore by the criterion established earlier, $\bar f$ is well-defined on $B(\eps_1(f))$.

Now we use a similar method to refine the choice of $\eps_1$ in order to achieve the second property in the statement of the lemma. For $f \in \Omega_1$, define $g_f: B(\eps_1(f)) \rightarrow [\frac{1}{8},+\infty)$ by
$$g_f = \left[\frac{1}{64}+ (\partial_{y_1}\bar{f}^1-1)^2 + (\partial_{y_2}\bar{f}^1)^2 + (\partial_{y_1}\bar{f}^2)^2 + (\partial_{y_2}\bar{f}^2-1)^2\right]^{\frac{1}{2}}.$$
It is easy to check that $\|d\bar{f}|_{y}-\id\| < g_f(y)$ for $y \in B(\eps_1(f))$, and it follows from the definition of $\Omega_1$ that $g_f(0) = \frac{1}{8}$.

Let $\gamma: (0,1] \times B(1) \rightarrow [0,1]$ be a smooth function such that $\gamma(\eps,y)$ is $1$ for $y \in B(\frac{1}{2}\eps)$, and $0$ for $y \notin B(\frac{3}{4}\eps)$. Then it is straightforward to show that the function $\gamma(\eps_1(f), \cdot)g_f$ is well-defined on all of $B(1)$, and it is smooth with respect to $f$. Recall that $c$ is the constant that appears in the Sobolev inequality (\ref{sobolev}). Define $$\eps(f) = \frac{1}{8c}\left[\frac{1}{16c^2\eps_1(f)^2} + \sum_{|\alpha| \leq 3} \int_{B(1)} |\partial^{\alpha}(\gamma(\eps_1(f), \cdot)g_f)|^2\right]^{-\frac{1}{2}}.$$

We claim that the function $\eps$ satisfies all the requirements in the statement of the lemma. First of all, $\eps(f)$ is well-defined since $\eps_1(f) > 0$. Also it is easy to see that $\eps(f) > 0$ and $\eps(f) \leq \frac{4c\eps_1(f)}{8c} = \frac{\eps_1(f)}{2} \leq \frac{1}{2}$. It follows that $\bar f$ is well-defined on $B(2\eps(f)) \subset B(\eps_1(f))$. Proving smoothness of $\eps$ is similar to proving smoothness of $\eps_1$. For $y \in B(\eps(f)) \subset B(\frac{1}{2}\eps_1(f))$, we have $\gamma(\eps_1(f),y)=1$. Then it follows from the mean value theorem that $$|g_f(y) - g_f(0)| \leq \sup_{B(1)}|d(\gamma(\eps_1(f), \cdot)g_f)|\cdot|y|.$$  Therefore for $y$ such that $|y| < \eps(f)$,
\begin{align*}
\|d\bar f|_{y}-\id\|
&< g_f(y) \\
&\leq g_f(0) + \sup_{B(1)} |d(\gamma(\eps_1(f), \cdot)g_f)| \cdot |y| \\
&\leq \frac{1}{8} + \frac{\sup_{B(1)} |d(\gamma(\eps_1(f), \cdot)g_f)|}{8c\left[\frac{1}{16c^2\eps_1(f)^2} + \sum_{|\alpha| \leq 3} \int_{B(1)} |\partial^{\alpha}(\gamma(\eps_1(f), \cdot)g_f)|^2\right]^{\frac{1}{2}}} \\
&\leq \frac{1}{8} + \frac{\|\gamma(\eps_1(f), \cdot)g_f\|_{C^1\left(\overline{B(1)}\right)}}{8c\|\gamma(\eps_1(f), \cdot)g_f\|_{H^3(B(1))}} \leq \frac{1}{8} + \frac{1}{8} = \frac{1}{4}.
\end{align*}
This completes the proof of Lemma \ref{epsilon}.
\end{proof}

Recall that $p: \SSS^2 \sm \{-x_0\} \to \RR^2$ is the stereographic projection from the North Pole, and $\bar f = p \circ f \circ p^{-1}$ is the local expression of the map $f$. We homotope each $f \in \Omega_1$ to a diffeomorphism whose restriction to $p^{-1}\left(B(\frac{1}{2}\eps(f))\right)$ is the identity map.

\begin{lemma}\label{Omega2}
There exists a smooth homotopy $S:I\times\Omega_1 \to \Omega_1$ such that for each $(t,f) \in I\times\Omega_1$,
\begin{enumerate}
\item
$S_0(f) = f$,
\item
$S_1(f)$ restricted to the neighborhood $p^{-1}\left(B(\frac{1}{2}\eps(f))\right)$ is the identity map,
\item
$S_t(\id_{\SSS^2}) = \id_{\SSS^2}$.
\end{enumerate}
\end{lemma}

\begin{proof}

Let $\gamma: (0, \frac{1}{2}] \times B(1) \rightarrow [0,1]$ be a smooth function such that $\gamma(\eps,y)$ is $1$ for $y \in B(\frac{1}{2}\eps)$ and $0$ for $y \notin B(\eps)$, and additionally $|\partial_y \gamma| < \frac{3}{\eps}$ everywhere.

Given $t \in [0,1]$ and $f \in \Omega_1$, we define $S_t(f)$ by cases.
\begin{itemize}
\item
For $y \in B(2\eps(f))$, the local expression $\overline{S_t(f)}$ is defined by $$\overline{S_t(f)}(y) = (1-t) \bar{f}(y) + t\left[\gamma(\eps(f), y)y + (1-\gamma(\eps(f), y))\bar{f}(y)\right].$$
\item
For $x \notin p^{-1}(B(\eps(f)))$, define $S_t(f)(x)=f(x)$.
\end{itemize}

We first check that each $S_t(f)$ is a well-defined smooth map. First of all, for $y$ such that $|y| < 2\eps(f)$, it follows from Lemma \ref{epsilon} that $\bar f(y)$ is well-defined. Thus $\overline{S_t(f)}$ is well-defined in the first case. The overlap of the two cases is when $\eps(f) \leq |y| < 2\eps(f)$. For such $y$, it follows from the choice of $\gamma$ that $\gamma(\eps(f),y)=0$ and $\overline{S_t(f)}(y) = \bar f(y)$. This shows that the definitions agree on the overlap. In each case, it is clear that $S_t(f)$ is smooth. The agreement on the overlap implies that $S_t(f)$ is a well-defined smooth map.

Now we show that each $S_t(f)$ is an orientation-preserving local diffeomorphism. It is enough to prove that the matrix $d \overline{S_t(f)}|_{y}$ is invertible and of positive determinant when $y \in B(\eps(f))$. For such $y$, an easy computation shows that
$$d \overline{S_t(f)}|_{y} = \id + \left[1-t\gamma(\eps(f),y)\right](d \bar{f}|_{y} - \id) + t(y-\bar{f}(y)) \cdot \partial_y \gamma(\eps(f),y).$$
To show that the matrix $d \overline{S_t(f)}|_{y}$ is invertible, it suffices to show that the difference $d \overline{S_t(f)}|_{y} - \id$ has operator norm less than $1$. The mean value theorem and Lemma \ref{epsilon} imply that
$$|\bar{f}(y)-y| \leq \sup_{B(\eps(f))}\|d \bar{f} - \id\|\cdot|y| < \frac{\eps(f)}{4}.$$ Then it follows from $|\partial_y \gamma| < \frac{3}{\eps}$ that
$$\| \left[1-t\gamma(\eps(f),y)\right](d \bar{f}|_{y} - \id) + t(y-\bar{f}(y)) \cdot \partial_y \gamma(\eps(f),y) \| < \frac{1}{4} + \frac{\epsilon(f)}{4} \frac{3}{\epsilon(f)} = 1.$$
Hence $d \overline{S_t(f)}|_{y}$ is invertible. Moreover, since $f \in \Omega_1$, the matrix $d \overline{S_0(f)}|_{y} = d\bar f|_{y}$ has positive determinant. Thus each $d \overline{S_t(f)}|_{y}$ has positive determinant since it depends continuously on $t$.

We can use the following standard topological argument to show that $S_t(f)$ is in fact a diffeomorphism. By using the compactness of $\SSS^2$ and the fact that $S_t(f)$ is a local diffeomorphism, we can conclude that $S_t(f)$ is a covering map. It follows from Theorem 5.1 on page 147 of \cite{bredon} that $S_t(f)$ is injective. Therefore $S_t(f)$ is a diffeomorphism. It is easy to check that $\overline{S_t(f)}(0) = 0$ and $d\overline{S_t(f)}|_{0} = \id$. Thus $S_t(f) \in \Omega_1$.

Furthermore, it follows from the smoothness of the function $\eps$ that $S_t(f)$ is smooth with respect to $(t,f)$. It is easy to see that $S_0(f) = f$. For each $y \in B(\frac{1}{2}\eps(f))$, it follows from $\gamma(\eps(f),y)=1$ that $\overline{S_t(f)}(y) = y$. Lastly, it is clear that $S_t(\id_{\SSS^2}) = \id_{\SSS^2}$.
\end{proof}

In Lemma \ref{Omega2}, we homotope each $f \in \Omega_1$ to the diffeomorphism $S_1(f)$ which is the identity map on the neighborhood $p^{-1}\left(B(\frac{1}{2}\eps(f))\right)$ of the South Pole. To complete the proof of Theorem \ref{Id}, we need to homotope it to the identity map on the complement of this neighborhood.

\begin{proof}[Proof of Theorem \ref{Id}]
Let $\tilde p: \SSS^2 \sm \{x_0\} \to \RR^2$ be the stereographic projection from the South Pole. We use $\tilde p$ as the coordinate chart throughout this proof. For each $f \in \Omega_1$, we have $S_1(f)(x_0) = x_0$. Thus the local expression
$$\overline{S_1(f)} = \tilde p \circ S_1(f) \circ \tilde p^{-1}$$
is a well-defined diffeomorphism of $\RR^2$. It follows from Lemma \ref{Omega2} that $\overline{S_1(f)}$ restricted to the open set $\tilde p \circ p^{-1}\left(B(\frac{1}{2}\eps(f))\right)$ is the identity map. It is clear that
$$\tilde p \circ p^{-1}\left(B\left(\frac{1}{2}\eps(f)\right)\right) = \{ y\in \RR^2: |y| > \frac{2}{\eps(f)} \}.$$
Let $\Psi_f: \left[-\frac{3}{\eps(f)},\frac{3}{\eps(f)}\right]^2 \to [-1,1]^2$ be the scaling map, given by $y \mapsto \frac{\eps(f)}{3}y$. Then it is easy to see that $\Psi_f \circ \overline{S_1(f)} \circ \Psi_f^{-1}$ is an orientation preserving diffeomorphism of $[-1,1]^2$, which is the identity map on a neighborhood of the boundary. Thus it belongs to the set $\mathcal F$ in Proposition \ref{F}. Recall that $F: I \times \mathcal F \to \mathcal F$ is the smooth strong deformation retraction to $\{\id_{[-1,1]^2}\}$. Define the homotopy $T$ in terms of the local expression as follows.
$$\overline{T_t(S_1(f))}(y)=
\begin{cases}
\Psi_f^{-1} \circ F_t(\Psi_f \circ \overline{S_1(f)} \circ \Psi_f^{-1}) \circ \Psi_f(y) & \text{if } y \in \left[-\frac{3}{\eps(f)},\frac{3}{\eps(f)}\right]^2, \\
y & \text{otherwise}.
\end{cases}
$$

It is clear that each $T_t(S_1(f))$ is an orientation preserving diffeomorphism of $\SSS^2$. Moreover, properties of the homotopy $F$ imply that $T_0(S_1(f)) = S_1(f)$ and $T_1(S_1(f)) = \id_{\SSS^2}$. The smoothness property of maps $F$, $S$, and $\eps$ implies that $T_t(S_1(f))$ is smooth with respect to $(t,f)$. Lastly, $T_t(\id_{\SSS^2}) = \id_{\SSS^2}$ since $F_t(\id_{[-1,1]^2}) = \id_{[-1,1]^2}$.

To finish the proof, we need to smoothly concatenate the homotopies $S$ and $T$. Let $\beta_1:[0,\frac{1}{2}] \to [0,1]$ be a smooth map which is $0$ in a neighborhood of $0$ and $1$ in a neighborhood of $\frac{1}{2}$. Similarly, let $\beta_2:[\frac{1}{2},1] \to [0,1]$ be a smooth map which is $0$ in a neighborhood of $\frac{1}{2}$ and $1$ in a neighborhood of $1$. Define the homotopy $R: I \times \Omega_1 \to \Omega_1$ by
$$R_t(f) =
\begin{cases}
S_{\beta_1(t)}(f) & t \in [0,\frac{1}{2}], \\
T_{\beta_2(t)}(S_1(f)) & t \in [\frac{1}{2},1].
\end{cases}$$
It is easy to see that $R$ satisfies all the required properties.
\end{proof}

To finish the proof of the main theorem, we need to smoothly concatenate the homotopy $Q$ of Theorem \ref{Omega1} and the homotopy $R$ of Theorem \ref{Id}.
\begin{proof}[Proof of Theorem \ref{main}]
Let $\beta_1:[0,\frac{1}{2}] \to [0,1]$ be a smooth map which is $0$ in a neighborhood of $0$ and $1$ in a neighborhood of $\frac{1}{2}$. Similarly, let $\beta_2:[\frac{1}{2},1] \to [0,1]$ be a smooth map which is $0$ in a neighborhood of $\frac{1}{2}$ and $1$ in a neighborhood of $1$. Define the homotopy $P: I \times \Omega \to \Omega$ by
$$P_t(f) =
\begin{cases}
Q_{\beta_1(t)}(f) & t \in [0,\frac{1}{2}], \\
\alpha(f)^{-1} \circ R_{\beta_2(t)}(\alpha(f) \circ Q_1(f)) & t \in [\frac{1}{2},1].
\end{cases}$$
It is easy to see that $R$ satisfies all the required properties.
\end{proof}

\

\section{Diffeomorphisms of the Square}

In this section we prove Proposition \ref{F}. First of all, we observe that $[-1,1]^2$ is diffeomorphic to $I^2$ as manifolds with corners. It is more convenient to work with $I^2$ than $[-1,1]^2$. Therefore by abuse of notation, we denote by $\mathcal F$ \emph{the space of those orientation-preserving diffeomorphisms of the square $I^2$ such that for each $f \in \mathcal F$, there exists a neighborhood of the boundary $\partial I^2$ on which $f$ is the identity map.}

To prove Proposition \ref{F}, it is equivalent to proving
\begin{theorem} \label{square}
There is a smooth strong deformation retraction $F: I \times \mathcal F \to \mathcal F$ to $\{ \id_{I^2} \}$. More precisely, for each $(t,f) \in I \times \mathcal F$,
\begin{enumerate}
\item
$F_0(f) = f$,
\item
$F_1(f) = \id_{I^2}$,
\item
$F_t(\id_{I^2}) = \id_{I^2}$.
\end{enumerate}
\end{theorem}

To construct the desired deformation retraction $F$, we consider a superset $\mathcal E$ of $\mathcal F$, and construct a deformation retraction from this superset to $\{ \id_{I^2} \}$. The set $\mathcal E$ is defined as follows.

We denote by $\e_1$ the vector $(1,0)$. Let $I_1$ be the right boundary of $I^2$. In other words, $I_1 = \{1\} \times I$. We denote by $\mathcal E$ \emph{the space of those orientation-preserving diffeomorphisms of the square $I^2$ such that for each $f \in \mathcal E$,}
\begin{itemize}
\item
\emph{there exists a neighborhood of $\partial I^2 \sm I_1$ on which $f$ is the identity map, and}
\item
\emph{for $x$ close enough to $1$, $df|_{(x,y)}\e_1 = \e_1$.}
\end{itemize}

By abuse of notation, we also let $\e_1$ denote the constant map taking $I^2$ to the vector $\e_1$. We view $\e_1$ as a constant vector field on $I^2$. Then for each $f \in \mathcal E$, there is a corresponding non-vanishing vector field $f_*\e_1$, i.e.,
$$(f_*\e_1)(x,y) = df|_{f^{-1}(x,y)}\e_1.$$
To homotope an element $f \in \mathcal E$ to the identity map, we first homotope its corresponding vector field $f_*\e_1$ to the constant vector field $\e_1$ (Lemma \ref{vector}), and then integrate the vector fields to recover the corresponding diffeomorphisms in $\mathcal E$ (Theorem \ref{superset}). We start by considering the following space of vector fields on $I^2$.

We denote by $\mathcal H$ \emph{the space of all $C^\infty$ maps from $I^2$ to $\RR^2 \sm \{0\}$ such that for each $h \in \mathcal H$, there exists a neighborhood of $\partial I^2$ on which $h$ is equal to the constant map $\e_1$.}

\begin{lemma} \label{vector}
There is a smooth homotopy $\Phi: I \times \mathcal E \to \mathcal H$, such that for each $(t,f) \in I \times \mathcal E$,
\begin{enumerate}
\item
$\Phi_0(f) = f_*\e_1$,
\item
$\Phi_1(f) = \e_1$,
\item
$\Phi_t(\id_{I^2}) = \e_1$.
\end{enumerate}
\end{lemma}

\begin{proof}
The map $\exp: \CC \to \CC \sm \{0\}$ is a $C^\infty$ covering map of $\RR^2 \sm \{0\}$. Fix a vector $\tilde \e_1 \in \exp^{-1}(\{ \e_1 \})$ in $\RR^2$. By Theorem 4.1 on page 143 of \cite{bredon}, for each $h \in \mathcal H$, there is a unique continuous map $\tilde h:I^2 \to \RR^2$ such that
$$\exp \circ \tilde h = h \text{ and } \tilde h(0,0) = \tilde \e_1.$$
Moreover, denote by $U_h$ a connected neighborhood of $\partial I^2$ on which $h$ is equal to $\e_1$. By going through the construction of $\tilde h$ in Theorem 4.1 of \cite{bredon}, one can show that $\tilde h (U_h) = \{ \tilde \e_1 \}$. Lastly, for a small enough neighborhood $V \subset I^2$, the image $h(V)$ is contained in a basic open set $W$, on which $\log$ is defined as a multi-valued map. Then there is a unique branch of $\log$ such that $\tilde h|_V = \log \circ h|_V$. Thus the map $\tilde h$ is $C^\infty$. A similar argument shows that the map $h \mapsto \tilde h$ is smooth.

Let $H: I \times \RR^2 \to \RR^2$ be a $C^\infty$ homotopy that contracts $\RR^2$ to $\{\tilde \e_1\}$. Then define $\Phi: I \times \mathcal E \to \mathcal H$ as follows. For $(t,f) \in I \times \mathcal E$,
$$\Phi_t(f) = \exp \circ H_t \circ \widetilde{f_*\e_1}.$$
By using the previously established properties, one can show that each $\Phi_t(f)$ is indeed an element of $\mathcal H$, the map $\Phi$ is smooth, and conditions (1) - (3) are satisfied.
\end{proof}

Now we show that the space $\mathcal E$ can be smoothly deformed to $\{ \id_{I^2} \}$.
\begin{theorem} \label{superset}
There is a smooth strong deformation retraction $E: I \times \mathcal E \to \mathcal E$ to $\{ \id_{I^2} \}$. More precisely, for each $(t,f) \in I \times \mathcal E$,
\begin{enumerate}
\item
$E_0(f) = f$,
\item
$E_1(f) = \id_{I^2}$,
\item
$E_t(\id_{I^2}) = \id_{I^2}$.
\end{enumerate}
\end{theorem}

\begin{proof}
For $(t,f) \in I \times \mathcal E$ and $y \in I$, let $s \mapsto \Gamma_t(f)(s,y)$ be the integral curve of the vector field $\Phi_t(f)$ from Lemma \ref{vector}, with the initial condition $\Gamma_t(f)(0,y) = (0,y)$. It follows from the standard regularity theory of ordinary differential equations that $\Gamma_t(f)(s,y)$ depends smoothly on $(t,f,s,y)$.

Firstly, for any $t \in I$, it follows from $\Phi_t(\id_{I^2}) = \e_1$ that $\Gamma_t(\id_{I^2})(s,y) = (s,y)$. Hence $\Gamma_t(\id_{I^2}) = \id_{I^2}$. When $t = 1$, the fact that $\Phi_1(f) = \e_1$ implies that $\Gamma_1(f) = \id_{I^2}$. When $t = 0$, by using the fact that $\Phi_0(f) = f_*\e_1$, one can check directly that $\Gamma_0(f)(s,y) = f(s,y)$. Thus $\Gamma_t(f)$ satisfies conditions (1) - (3). However, $(t,f) \mapsto \Gamma_t(f)$ is not the desired homotopy; for $0<t<1$, we need to analyze $\Gamma_t(f)$ carefully.

Since each $\Phi_t(f)$ agrees with $\e_1$ on a neighborhood of $\partial I^2$, we can conclude that the integral curve $\Gamma_t(f)(\cdot,y)$ either
\begin{itemize}
\item
leaves the square $I^2$ via the right boundary $I_1$, or
\item
does not leave $I^2$.
\end{itemize}
However, in the second case, the integral curve would approach asymptotically to a simple closed curve, then by the Poincar\'e-Bendixson theorem (see page 191 of \cite{robinson}), the vector field $\Phi_t(f)$ vanishes somewhere in the interior of the closed curve, which contradicts the assumption that $\Phi_t(f) \in \mathcal H$. Thus each integral curve $\Gamma_t(f)(\cdot,y)$ meets $I_1$ at some time $\bar s$, which depends on $(t,f,y)$. More precisely, $\bar s(t,f,y)$ is defined by the equation
$$\Gamma_t(f)^1(\bar s,y) = 1.$$
To show that $\bar s$ is smooth, we compute the partial derivative
\begin{equation} \label{partial}
\frac{\partial}{\partial s}\Gamma_t(f)^1(s,y)\Big|_{s=\bar s} = \Phi_t(f)^1\left(\Gamma_t(f)(\bar s,y)\right) = 1,
\end{equation}
where the last equality follows from $\Phi_t(f)|_{I_1} = \e_1$. Therefore it follows from the implicit function theorem that $\bar s(t,f,y)$ depends smoothly on $(t,f,y)$.

Note that $\Gamma_t(f)$ is a diffeomorphism from the set $\bigcup_{y} \{ y \} \times \left[0,\bar s(t,f,y)\right]$ to the square $I^2$. However $\bar s$ is not necessarily equal to 1, so we cannot take $\Gamma_t(f)$ as the desired homotopy. We need to reparametrize the $s$ variable in $\Gamma_t(f)(s,y)$ to obtain an element of $\mathcal E$.

Let $\chi: [0,1] \times (0,+\infty) \to [0,+\infty)$ be a smooth function, such that for each $r \in (0,+\infty)$,
\begin{itemize}
\item
the function $\chi(\cdot,r)$ maps the interval $[0,1]$ diffeomorphically to $[0,r]$, with $\chi(0,r)=0$ and $\chi(1,r)=r$,
\item
there exists a neighborhood $U_0$ of $0$ and a neighborhood $U_1$ of $1$, such that for $x \in U_0 \cup U_1$, we have $\frac{\partial}{\partial x}\chi(x,r)=1$,
\item
for $x \in [0,1]$, we have $\chi(x,1)=x$.
\end{itemize}

We use the diffeomorphism $x \mapsto s=\chi(x,\bar s)$ to reparametrize $\Gamma_t(f)(s,y)$. We claim that
$$E_t(f)(x,y) := \Gamma_t(f)\left(\chi\left(x,\bar s(t,f,y)\right),y\right)$$
is the right deformation retraction.

We first show that each map $E_t(f):I^2 \to I^2$ has a smooth inverse. Fix $(t,f) \in I \times \mathcal E$. Let $s \mapsto \gamma_t(f)(s,x',y')$ be the integral curve of the vector field $\Phi_t(f)$, with the initial condition $\gamma_t(f)(0,x',y')=(x',y')$. Let us denote by $-\tau(x',y')$ the time when the integral curve $\gamma_t(f)(\cdot,x',y')$ meets the left boundary. More precisely, $-\tau(x',y')$ is defined by the equation
$$\gamma_t(f)^1(-\tau,x',y')=0.$$
Define $y(x',y')$ by $y=\gamma_t(f)^2(-\tau,x',y')$. It is easy to see that $\Gamma_t(f)(\tau,y)=(x',y')$. Now it suffices to find $x(x',y')$. Recall that by the first property of $\chi$, the function $\chi(\cdot,\bar s(y)):[0,1] \to [0,\bar s(y)]$ has a smooth inverse function, which we denote by $\sigma(y)$. Define $x(x',y')$ by $x=\sigma(y)(\tau)$. Then it is easy to check that the inverse of $E_t(f)$ is given by
$$(E_t(f))^{-1}(x',y') = \left(x(x',y'),y(x',y')\right).$$
To show that $(E_t(f))^{-1}$ is smooth, it suffices to show the function $\tau$ is smooth. One can carry out a similar computation as equation (\ref{partial}), and then the smoothness of $\tau$ follows from the implicit function theorem.

One can check that each $E_t(f)$ is in $\mathcal E$ by using the second property of the function $\chi$ and the property of $\Phi_t(f)$ near the boundary. Lastly, the fact that $\Gamma_t(f)$ satisfies conditions (1) - (3) and the third property of the function $\chi$ implies that $E: I \times \mathcal E \to \mathcal E$ is the desired deformation retraction. This completes the proof of Theorem \ref{superset}.
\end{proof}

\begin{remark} \label{comparison3}
On the top of page 623 of \cite{smale}, Smale uses an argument in general topology to reparametrize the variable $s$, which yields a continuous but not necessarily smooth homotopy. Here we use the smooth function $\chi$ to reparametrize, and the resulting homotopy is smooth.
\end{remark}

To finish the construction of the homotopy from $\mathcal F$ to $\{ \id_{I^2} \}$, we first notice that we \emph{cannot} define the homotopy $F: I \times \mathcal F \to \mathcal F$ to be the restriction of the homotopy $E$. This is because for $f \in \mathcal F$, the diffeomorphism $E_t(f)$ does not necessarily lie in the set $\mathcal F$. To solve this problem, we first construct a retraction $p:\mathcal E \to \mathcal F$, and we show that $F_t(f)$ given by $(p \circ E_t)(f)$ is the right homotopy. To construct this retraction, we first consider the following set of diffeomorphisms of the interval $I$.

Let $\mathcal G$ be \emph{the space of those orientation-preserving diffeomorphisms of $I$ such that for each $g \in \mathcal G$, there exists a neighborhood of $0$ and a neighborhood of $1$ on which $g$ is the identity map.}

\begin{lemma} \label{interval}
There is a smooth strong deformation retraction $G: I \times \mathcal G \to \mathcal G$ to $\{ \id_I \}$. More precisely, for each $(t,g) \in I \times \mathcal G$,
\begin{enumerate}
\item
$G_0(g) = g$,
\item
$G_1(g) = \id_I$,
\item
$G_t(\id_I) = \id_I$.
\end{enumerate}
\end{lemma}

\begin{proof}
Define the homotopy $G: I \times \mathcal G \to \mathcal G$ as follows. For $x \in I$,
$$G_t(g)(x) = (1-t)g(x) + tx.$$
It is easy to check that $G$ is smooth and conditions (1) - (3) are satisfied.
\end{proof}

\begin{proof}[Proof of Theorem \ref{square}]
For each $f \in \mathcal E$, one can check that for $x$ close enough to $1$, the value $f(x,y)$ is given by $(x,g_f(y))$, where $g_f(y)=f^2(1,y)$. It is easy to see that $g_f$ is an element of $\mathcal G$ and it depends smoothly on $f$.

Let $\beta: I \to I$ be a smooth function that is 1 in a neighborhood of 0, and 0 in a neighborhood of 1. Now define $\Psi_f: I^2 \to I^2$ by
$$\Psi_f(x,y) = \left(x, \left[G_{\beta(x)}\left(g_f^{-1}\right)\right](y)\right).$$
It is easy to check that the inverse of $\Psi_f$ is given by the map $(x,y) \mapsto \left(x, \left[G_{\beta(x)}\left(g_f^{-1}\right)\right]^{-1}(y)\right)$. Thus $\Psi_f$ is a diffeomorphism of the square.

One can show that the map $p: \mathcal E \to \mathcal F$  defined by
$$p(f) = \Psi_f \circ f$$
is a smooth retraction. In other words, $p$ is smooth and $p|_{\mathcal F} = \id$. It follows that the homotopy $F: I \times \mathcal F \to \mathcal F$ defined by
$$F_t(f) = (p \circ E_t)(f)$$
is the desired smooth deformation retraction to $\{ \id_{I^2} \}$.
\end{proof}

\


\begin{thebibliography}{9}

\bibitem{bredon}
G. Bredon, \emph{Topology and Geometry}, Graduate Texts in Mathematics, vol 139, 1993

\bibitem{earle}
C.J. Earle, J. Eells, The Diffeomorphism Group of a Compact Riemann Surface, \emph{Bull. Amer. Math. Soc.} 73 (1967), 557-559

\bibitem{evans}
L. Evans, \emph{Partial Differential Equations}, Graduate studies in mathematics, vol 19, 1998

\bibitem{iglesias}
P. Iglesias-Zemmour, \emph{Diffeology}, http://math.huji.ac.il/$\sim$piz/documents/Diffeology.pdf

\bibitem{joyce}
D. Joyce, On manifolds with corners, http://arxiv.org/abs/0910.3518

\bibitem{robinson}
R.C. Robinson, \emph{An Introduction to Dynamical Systems}, Person Prentice Hall, New Jersey, 2004.

\bibitem{schwarz}
G. Schwarz, Smooth functions invariant under the action of a compact Lie group, \emph{Topology} 14 (1975), 63-68

\bibitem{smale}
S. Smale, Diffeomorphisms of the 2-Sphere, \emph{Proceedings of the American Mathematical Society} 10 (1959), no. 4, 621-626

\end{thebibliography}
\end{document}